\documentclass[11pt]{article}
\usepackage{amscd}
\usepackage{amsfonts}
\usepackage{amsmath}
\usepackage{amssymb}
\usepackage{amsthm}
\usepackage{bbm}
\usepackage{CJK}
\usepackage{fancyhdr}
\usepackage{graphicx}
\usepackage{hyperref}
\usepackage{indentfirst}
\usepackage{latexsym}
\usepackage{mathrsfs}
\usepackage{xypic}

\newtheorem{theorem}{Theorem}[section]
\newtheorem{lemma}[theorem]{Lemma}
\newtheorem{definition}[theorem]{Definition}
\newtheorem{proposition}[theorem]{Proposition}
\newtheorem{example}[theorem]{Example}

\newtheorem{cor}[theorem]{Corollary}

\usepackage[top=1in,bottom=1in,left=1.25in,right=1.25in]{geometry}
\def\<{\langle}
\def\>{\rangle}

\def\c{\cdot}

\def\o{\otimes}

\date{}
\begin{document}
\renewcommand{\baselinestretch}{1.2}
\renewcommand{\arraystretch}{1.0}
\title{\bf  Rota-Baxter operators and  related structures on anti-flexible algebras}
\date{}
\author{{\bf Shuangjian Guo$^{1}$, Ripan Saha$^{2}$\footnote
        { Corresponding author:~~ripanjumaths@gmail.com} }\\
{\small 1. School of Mathematics and Statistics, Guizhou University of Finance and Economics} \\
{\small  Guiyang  550025, P. R. of China} \\
{\small 2. Department of Mathematics, Raiganj University } \\
{\small  Raiganj, 733134, West Bengal, India}}
 \maketitle
\begin{center}
\begin{minipage}{13.cm}

{\bf \begin{center} ABSTRACT \end{center}}
In this paper, we first construct a graded Lie algebra which characterizes  Rota-Baxter operators on an anti-flexible algebra as Maurer-Cartan elements. Next, we study infinitesimal deformations
of bimodules over anti-flexible algebras. We also consider compatible Rota-Baxter operators on bimodules over anti-flexible algebras. Finally, We define $\mathcal{ON}$-structures which give rise to compatible Rota-Baxter operators and vice-versa.
 \smallskip

 \smallskip

 \end{minipage}
 \end{center}
 \normalsize\vskip0.5cm

\section{Introduction}
\def\theequation{\arabic{section}. \arabic{equation}}
\setcounter{equation} {0}

Flexible algebra was introduced by Oehmke \cite{O58} as a natural generalization of associative algebras. Let $A$ be a vector space over a field $\mathbb{K}$ equipped with a bilinear product $A\times A\to A,~~~(a, b) \mapsto ab$. We denote the associator of $A$ as
\begin{eqnarray*}
 (a, b, c) = (ab)c - a(bc), ~\text{for all}~a, b, c \in A.
\end{eqnarray*}
$A$ is called a flexible algebra if the following identity is satisfied
\begin{eqnarray*}
(a, b, a) = 0, \mbox{or equivalently}, (ab)a = a(ba), ~\text{for all}~a, b \in A.
\end{eqnarray*}
Anti-flexible algebras are a natural generalization of flexible algebras introduced by Rodabaugh in \cite{R65}. Rodabaugh studied anti-flexible algebras in more detail in \cite{R67, R72, R78}. Let $A$ be a vector space equipped with a bilinear product $(a, b) \mapsto a\c b$. $A$ is called an anti-flexible algebra if the following identity is satisfied
\begin{eqnarray*}
(a, b, c) = (c, b, a), \mbox{or equivalently},   (a\c b)\c c- a\c (b\c c) = (c\c b)\c a- c\c (b\c a), ~\text{for all}~ a, b, c \in A.
\end{eqnarray*}
The notion of anti-flexible bialgebras was studied in \cite {DBH20}. Goze and Remm \cite{GR04} constructed a graded Lie algebra structure on the graded space of all multilinear maps on a vector space, and studied cohomology and the deformation of anti-flexible algebras. In this paper,  we construct a graded Lie algebra structure that characterizes Rota-Baxter operators on an anti-flexible algebra as Maurer-Cartan elements.

The deformation of algebraic structures began with the seminal work of Gerstenhaber \cite{G63, G64} for associative algebras, and followed by its extension to Lie algebras by Nijenhuis and Richardson \cite{NR66, NR68}.  In general, deformation theory was developed for binary quadratic operads by Balavoine \cite{B97}. Deformations of morphisms were
developed in \cite{FZ15, TBGS19}.  

While studying the fluctuation theory in probability, the notion of Rota-Baxter operators on associative algebras was introduced by Baxter \cite{B60} in
1960. Since its introduction, it has been found many applications, including in Connes-Kreimer's algebraic approach
to the renormalization in perturbative quantum field theory \cite{CK00}. Rota-Baxter operator is closely related with dendriform algebras, Lie algebras, and solution of the classical Yang-Baxter equation, see \cite{L01, G12} for more details. Rota-Baxter operators are also useful in the study of
dendriform algebras operads, which give rise to the splitting of operads \cite{BBGN13, PBG17}. With motivation from Poisson structures, the notion of Rota-Baxter operators on bimodules over associative algebras was introduced by Uchino \cite{U08}.  Recently, the notions of compatible Rota-Baxter operators and $\mathcal{ON}$-structures was introduced by Liu, Bai, and Sheng in \cite{LBS19}, and they proved that an $\mathcal{ON}$-structure gives rise to a hierarchy of Rota-Baxter operators, and that a solution of the strong Maurer-Cartan equation on the associative twilled algebra associated to a Rota-Baxter operator gives rise to a pair of $\mathcal{ON}$-structures which are naturally in duality. In \cite{Da20}, Das constructed an explicit graded Lie algebra whose Maurer-Cartan elements are Rota-Baxter operators on associative algebras and studied linear and formal deformations of a Rota-Baxter operator on an associative algebra. Our main objectives in this paper are certain operators on anti-flexible algebras. More precisely, we are interested in the notions of compatible Rota-Baxter operators and $\mathcal{ON}$-structures on anti-flexible algebras. We show that an $\mathcal{ON}$-structure gives rise to compatible Rota-Baxter operators and conversely given two compatible Rota-Baxter operators there is an $\mathcal{ON}$-structures such that this correspondence naturally in duality.

The paper is organized as follows. In Section \ref{sec 2}, we construct the graded Lie algebra that characterizes  Rota-Baxter operators on an anti-flexible algebra as Maurer-Cartan elements. In Section \ref{sec 3}, we show that the cohomology of a Rota-Baxter operator can also be described as the Hochschild cohomology of a certain anti-flexible algebra with a suitable bimodule. We also relate the cohomology of a Rota-Baxter operator on an anti-flexible algebra with the cohomology of the corresponding Rota-Baxter operator on the commutator Lie algebra. In Section \ref{sec 4},  we study infinitesimal deformations
of bimodules over anti-flexible algebras. In Section \ref{sec 5}, we consider compatible Rota-Baxter operators on bimodules over anti-flexible algebras. We define $\mathcal{ON}$-structures which give rise to a hierarchy of compatible Rota-Baxter operators.

Throughout this paper, we work over the complex field $\mathbb{K}$, and all the vector spaces are finite-dimensional.
\section{Rota-Baxter operators}\label{sec 2}
\def\theequation{\arabic{section}.\arabic{equation}}
\setcounter{equation} {0}

In this section, we first recall the basics of Rota-Baxter operators on anti-flexible algebra and their morphisms. Then, we construct a graded Lie algebra with a graded Lie bracket
whose Maurer-Cartan elements are Rota-Baxter operators on anti-flexible algebra. This construction allows us to define cohomology for a Rota-Baxter operator.
\begin{definition}(\cite{DBH20})
Let $A$ be a vector space equipped with a bilinear product $(x, y) \rightarrow x\c y$. $A$ is
called an anti-flexible algebra if the following identity is satisfied
\begin{eqnarray}
 (a\c b)\c c- a\c (b\c c) = (c\c b)\c a- c\c (b\c a), ~ \text{for all}~ a, b, c \in A.
\end{eqnarray}
\end{definition}

\begin{example}
Every associative algebra is automatically an anti-flexible algebra.
\end{example}
\begin{example}
Let $(A, \c)$ be an anti-flexible algebra and $B$ an associative algebra, then $A\o B$ is an  anti-flexible algebra with product given by
\begin{eqnarray*}
(a_1\o b_1) \c (a_2\o b_2)=a_1\c a_2\o b_1b_2,~ \text{for all}~  a_1,a_2\in A, b_1,b_2\in B.
\end{eqnarray*}
\end{example}
\begin{example}
Let $(A, \c)$ and $(B, \c)$  be  anti-flexible algebras, then $(A\oplus B, \c)$ is an anti-flexible algebra with the operation componentwise multiplication.
\end{example}

\begin{definition}(\cite{DBH20})
Let $(A, \c)$ be an anti-flexible algebra and $M$ be a vector space. Let $l, r : A \rightarrow
gl(M)$ be two linear maps. If for any $a, b \in A$,
\begin{eqnarray}
&&l(a\c b)-l(a)l(b) = r(a)r(b)-r(b \c a),\\
&&l(a)r(b)- r(b)l(a) = l(b)r(a)- r(a)l(b).
\end{eqnarray}
Then it is called a bimodule of $(A, \c)$, denoted by $(M, l, r)$.
\end{definition}

Given an anti-flexible algebra $(A, \c)$ and a bimodule $(M, l, r)$, the vector space $A\oplus M$ carries
an anti-flexible algebra structure with product given by
\begin{eqnarray*}
(a, m)\c (b, n)=(a\c b, l(a)n+r(b)m), ~ \text{for all}~ a, b\in A, m, n\in M.
\end{eqnarray*}
This is called the semi-direct product of $A$ with $M$.

\begin{definition}(\cite{DBH20})
Rota-Baxter operator on an anti-flexible algebra $(A, \c)$ with respect to the bimodule $(M, l, r)$ is given by a linear map $T : M \rightarrow A$
that satisfies
 \begin{eqnarray}
 T(m) \c T(n) = T(l(T(m))n + r(T(n))m), ~~~~~\forall m, n \in M.
 \end{eqnarray}
\end{definition}

Following Uchin\cite{U08}, we have the following proposition.
\begin{proposition}
A linear map $T : M \rightarrow  A$ is a Rota-Baxter operator on an anti-flexible algebra $A$ with respect to the bimodule $(M, l, r)$ if and only if the graph
\begin{eqnarray*}
Gr(T) = \{(T(m), m)| m \in M\}
\end{eqnarray*}
is a subalgebra of the semi-direct product algebra $A \oplus M$.
\end{proposition}

\begin{definition}
Let $(A, \c)$ be an anti-flexible algebra. A linear map $N : A \rightarrow A$ is said
to be a Nijenhuis operator if its Nijenhuis torsion vanishes, that is,
\begin{eqnarray*}
N(a) \c N(b) = N(Na \c b + a \c Nb - N(a \c b)),~ \text{for all}~ a, b \in A.
\end{eqnarray*}
The  operation $\c_N : A \o A\rightarrow  A$ given by
\begin{eqnarray*}
~~~a \c_N b = Na \c b + a \c Nb - N(a \c b),~ \text{for all}~ a, b \in A.
\end{eqnarray*}
is an anti-flexible algebra and $N$ is an anti-flexible algebra homomorphism from $(A, \c_N)$ to $(A, \c)$.
\end{definition}
By direct calculations, we have
\begin{lemma}
Let $(A, \c)$ be an anti-flexible algebra and $N$ be a Nijenhuis operator on A. For all $l, k \in  \mathbb{K}$,\\
(i) $(A, \c_{N^k} )$ is an anti-flexible algebra, \\
(ii) $N^l$ is also a Nijenhuis operator on the anti-flexible algebra $(A, \c_{N^k} )$,\\
(iii) The anti-flexible algebras $(A, (\c_{N^k})_{N^{l}} )$ and $(A, \c_{N^{k+l}})$ coincide,\\
(iv) The anti-flexible algebras $(A, \c_{N^{k}})$ and $(A, \c_{N^{l}})$ are compatible, that is, any linear combination of
$\c_{N^{k}}$ and $\c_{N^{l}}$ still makes $A$ into an anti-flexible algebra,\\
(v) $N^l$ is an anti-flexible algebra homomorphism from $(A, \c_{N^{k+l}})$ to $(A, \c_{N^{k}})$.
\end{lemma}

Another characterization of a Rota-Baxter operator can be given in terms of anti-flexible-Nijenhuis operator
on anti-flexible algebras.
\begin{proposition}
A linear map $T : M \rightarrow A$ is  a Rota-Baxter operator on $(A, \c)$ with respect to the
bimodule $(M, l, r)$ if and only if $N_T=\left(
                                           \begin{array}{cc}
                                             0 & T \\
                                             0 & 0 \\
                                           \end{array}
                                         \right): A \oplus M\rightarrow  A \oplus M$ is an anti-flexible-Nijenhuis operator on the semi-direct product algebra $A \oplus M$.
\end{proposition}

Next, we recall pre-anti-flexible algebra structures which were first introduced by Bai \cite{DBH20},  pre-anti-flexible algebras  can be regarded as a natural generalization of dendriform algebras introduced by Loday \cite{L01}. On the other hand, from
the point of view of operads, like dendriform algebras being the splitting of associative algebras, pre-anti-flexible algebras are the splitting of anti-flexible algebras (\cite{BBGN13,PBG17}).
\begin{definition}(\cite{DBH20})
Let $A$ be a vector space with two bilinear products $\prec, \succ: A \o A \rightarrow A$. We call it
a pre-anti-flexible algebra denoted by $(A, \prec, \succ)$ if for any $a, b, c \in A$, the following equations
are satisfied
\begin{eqnarray}
&& (a \succ b)\prec c - a \succ (b \prec c)= (c \succ b)\prec a - c \succ (b \prec a),\\
&&  (a \ast b) \succ c - a \succ (b \succ c)= (a \prec b)\prec c-a \prec (b \ast c),
\end{eqnarray}
where $a\ast b=a\prec b+a\succ b$.
\end{definition}
A Rota-Baxter operator has an underlying pre-anti-flexible algebra structure \cite{DBH20}.
\begin{proposition}
Let $T: M \rightarrow A$ be a Rota-Baxter operator on an anti-flexible algebra $(A, \c)$ with respect to the bimodule $(M, l, r)$.
Then the vector space $M$ carries a pre-anti-flexible algebra  structure with
\begin{eqnarray*}
m \succ n = l(T(m))n, ~~~~~~    m \prec n = r(T(n))m, ~ \text{for all}~  m, n \in M.
\end{eqnarray*}
\end{proposition}

\begin{definition}
A morphism of Rota-Baxter operator from $T$ to $T'$ consists of a pair $(\phi, \psi)$ of an
algebra morphism $\phi: A \rightarrow B$ and a linear map $\psi: M \rightarrow N$ satisfying
\begin{eqnarray}
T' \circ \psi = \phi \circ T,\\
l(\phi(a))\psi(m) = \psi(l(a)m), \\
r(\phi(a))\psi(m) = \psi(r(a)m),
\end{eqnarray}
for any $a\in A$ and $m\in M$.
\end{definition}
It is called an isomorphism if $\phi$ and $\psi$ are both linear isomorphisms.

The proof of the following result is straightforward and we omit the details.
\begin{proposition}
A pair of linear maps $(\phi: A \rightarrow B, ~\psi: M \rightarrow N)$ is a morphism of
Rota-Baxter operators from $T$ to $T'$ if and only if
\begin{eqnarray*}
Gr((\phi, \psi)):=\{((a, m),(\phi(a), \psi(m)))| a \in A, m \in M\}\subset (A \oplus M) \oplus (B \oplus N)
\end{eqnarray*}
is a subalgebra, where $A \oplus M$ and $B \oplus N$ are equipped with semi-direct product algebra structures.
\end{proposition}

\begin{proposition}
Let $T$ be a  Rota-Baxter operator on an anti-flexible algebra $(A, \c)$ with respect to a bimodule $(M, l, r)$ and $T'$ be
a Rota-Baxter operator on $(B, \c)$ with respect to a bimodule $(N, l, r)$. If $(\phi, \psi)$ is a morphism from $T$ to $T'$,
then $\psi : M \rightarrow N$ is a morphism between induced pre-anti-flexible algebra structures.
\end{proposition}
{\bf Proof.} For all $m, m' \in  M$, we have
\begin{eqnarray*}
\psi(m\prec_M m')&=& \psi(r(T(m'))m) \\
&=& r(\phi(T(m')))\psi(m)\\
&=& r(T'\phi(m'))\psi(m)\\
&=&\psi(m)\prec_N\psi( m').
\end{eqnarray*}
Similary, we obtain $\psi(m\succ_M m')=\psi(m)\succ_N\psi(m')$.  \hfill $\square$

In the sequel,  we follow the result of  Goze and Remm \cite{GR04}  and the derived bracket construction of Voronov \cite{V05}
to construct an explicit graded Lie algebra whose Maurer-Cartan elements are Rota-Baxter operators. This construction is somewhat similar to Das \cite{Da20} but more helpful to study deformation theory of Rota-Baxter operators.

Recall that, in \cite{GR04} Goze and Remm constructed a graded Lie algebra structure on the
graded space of all multilinear maps on a vector space $V$. Recall that, for each $n\geq 0$, $\mathbf{g}^n = \text{Hom}(V^{\otimes n+1}, V )$ and a graded Lie bracket on $\oplus_n\mathbf{g}^n$ by:
\begin{eqnarray*}
[f, g]=f\overline{\circ}g- (-1)^{mn}g\overline{\circ}f, ~\text{for all}~ f\in \mathbf{g}^{m}, g\in \mathbf{g}^{n},
\end{eqnarray*}
and $\overline{\circ}$ is defined by
\begin{eqnarray*}
&&(f\overline{\circ} g)(x_1,\ldots, x_{m+n+1})\\
&=&\sum_{i=1}^{m+1} \sum_{\sigma\in \Sigma_{m+n-1}}(-1)^{\epsilon (\sigma)}(-1)^{(i-1)(n-1)}f(x_{\sigma(1)}, \ldots, x_{\sigma(i-1)}, g(x_{\sigma(i)}, \ldots, x_{\sigma(i+n-1)}, x_{\sigma(i+n)}), \\
&&\hspace{10cm}x_{\sigma(i+m+1)}, \ldots, x_{\sigma(m+n+1)}),
\end{eqnarray*}
where $\Sigma_p$ refers to the $p$-symmetric group and $\epsilon(\sigma)$ denotes the sign of $\sigma$.

Let $A$ be an anti-flexible algebra equipped with multiplication map $\mu: A\o A\rightarrow A,~\mu(a,b)=a.b$. We know that $A\oplus M$ has also an anti-flexible algebra structure. Consider the graded Lie algebra structure on $\mathbf{g}^n = \text{Hom}((A\oplus M) ^{\otimes n+1}, A\oplus M )$ associated to the direct sum vector space $V = A \oplus M$. Observe that the elements $\mu, l, r\in \mathbf{g}^1 = \text{Hom}((A\oplus M) \otimes ^{2}, A\oplus M )$. Therefore, $\mu +l+r\in \mathbf{g}^1$.

\begin{proposition}
The product $\mu$ defines a multiplication structure on $A$ and $l, r$ defines an $A$-bimodule structure on $M$ if and only if $(\mu +l+r)\overline{\circ}(\mu +l+r)=0$, i.e. $\mu +l+r\in \mathbf{g}^1$ is a Maurer-Cartan element in  $\mathbf{g}$.
\end{proposition}
{\bf Proof.} For any $a_1, a_2, a_3\in A$ and $m_1, m_2, m_3\in M$,  we have
\begin{eqnarray*}
 &&(\mu +l+r)\overline{\circ}(\mu +l+r)((a_1, m_1),(a_2, m_2), (a_3, m_3))\\
 &=&  (\mu +l+r)((\mu +l+r)((a_1, m_1),(a_2, m_2)), (a_3, m_3))\\
 &&-(\mu +l+r)(((a_1, m_1),(\mu +l+r)((a_2, m_2)), (a_3, m_3)))\\
 &&-(\mu +l+r)((\mu +l+r)((a_3, m_3),(a_2, m_2)), (a_1, m_1))\\
 &&+(\mu +l+r)(((a_3, m_3),(\mu +l+r)((a_2, m_2)), (a_1, m_1)))\\
 &=& ((a_1a_2)a_3, l(a_1a_2)m_3+l(a_1)r(a_3)m_2+r(a_3)r(a_2)(m_1))\\
 && -(a_1(a_2a_3), l(a_1)l(a_2)m_3+r(a_3)l(a_1)m_2+ r(a_3)r(a_2)(m_1))\\
 &&-((a_3a_2)a_1, r(a_1)r(a_2)m_3+l(a_3)r(a_1)m_2+r(a_3)r(a_2)(m_1))\\
 && +(a_3(a_2a_1), r(a_2a_1)m_3+r(a_1)l(a_2)m_2+ r(a_3)r(a_2)(m_1))\\
 &=&0.
\end{eqnarray*}
This holds if and only if $\mu$ defines a multiplication structure on  the anti-flexible  algebra $A$ and $l, r$ define an
$A$-bimodule structure on $M$. \hfill $\square$

Consider the graded vector space
\begin{eqnarray*}
C^{\ast}(M, A):=\oplus_{n\geq 1}C^{n}(M, A)=\oplus_{n\geq 1} \text{Hom}(\otimes M^{\otimes n}, A).
\end{eqnarray*}

\begin{theorem}
With the above notations, $(C^{\ast}(M, A), [[\c, \c]])$ is a graded Lie algebra, where the
graded Lie bracket $[[\c, \c]]: C^{m}(M, A)\times C^{n}(M, A)\rightarrow C^{m+n}(M, A)$ is defined by
\begin{eqnarray*}
[[P,P']] := (-1)^{m}[[\mu+l+r,P],P'],
\end{eqnarray*}
for any $P\in C^{m}(M, A), P'\in C^{n}(M, A)$.

 More precisely, we have
\begin{eqnarray*}
&&[[P,P']](v_1, \ldots, v_{m+n})\\
&=& \sum_{k=1}^{m}\sum_{\sigma\in \Sigma_{m+n}}(-1)^{(k-1)n}(-1)^{\epsilon(\sigma)}P(v_{\sigma(1)}, \ldots, v_{\sigma(k-1)},\\
&&\hspace{6cm}l(P'(v_{\sigma(k)}, \ldots, v_{\sigma(k+n-1)})v_{\sigma(k+n)},\ldots, v_{\sigma(m+n)})\\
&&-\sum_{k=1}^{m}\sum_{\sigma\in \Sigma_{m+n}}(-1)^{kn}(-1)^{\epsilon(\sigma)}P(v_{\sigma(1)}, \ldots, v_{\sigma(k-1)},\\
&&\hspace{4cm}r(P'(v_{\sigma(k+1)}, \ldots, v_{\sigma(k+n-1)}))v_{\sigma(k)}, v_{\sigma(k+n+1)},\ldots, v_{\sigma(m+n)})\\
&&-\sum_{k=1}^{n}\sum_{\sigma\in \Sigma_{m+n}}(-1)^{(k+n-1)m}(-1)^{\epsilon(\sigma)}P'(v_{\sigma(1)}, \ldots, v_{\sigma(k-1)},\\
&&\hspace{5cm}l(P(v_{\sigma(k)}, \ldots, v_{\sigma(k+m-1)}))v_{\sigma(k+m)}, \ldots, v_{\sigma(m+n)})\\
&&+\sum_{k=1}^{n}\sum_{\sigma\in \Sigma_{m+n}}(-1)^{(k+n)m}(-1)^{\epsilon(\sigma)}P'(v_{\sigma(1)}, \ldots, v_{\sigma(k-1)},\\
&&\hspace{4cm}r(P(v_{\sigma(k+1)}, \ldots, v_{\sigma(k+m)}))v_{\sigma(k)}, v_{\sigma(k+m+1)}, \ldots, v_{\sigma(m+n)})
\end{eqnarray*}
\begin{eqnarray*}
&&+ (-1)^{mn}[P(v_{1}, \ldots, v_{m})P'(v_{m+1}, \ldots, v_{m+n})\\
&&\hspace{5cm}-(-1)^{mn}[P'(v_{1}, \ldots, v_{n-1},v)_{n})P(v_{n+1}, \ldots, v_{m+n})],
\end{eqnarray*}
for any $P\in C^{m}(M, A), P'\in C^{n}(M, A)$.
Moreover, its Maurer-Cartan elements are Rota-Baxter operator on the anti-flexible algebra
$(A, \c)$ with respect to the bimodule $(M, l, r)$.
\end{theorem}
{\bf Proof.} The graded Lie algebra  $(C^{\ast}(M, A), [[\c, \c]])$  is obtained via the derived bracket \cite{GR04}.   In fact,
the Balavoine bracket $[\c, \c]$ associated to the direct sum vector space $ A\oplus M$ gives rise to a graded
Lie algebra $(C^{\ast}(A\oplus M,  A\oplus M), [\c, \c])$. By
the above proposition,  we deduce that $(C^{\ast}(A\oplus M,  A\oplus M), [\c, \c], d=[\mu+l+r,\c])$  is a differential graded Lie algebra. Obviously $C^{\ast}(M, A)$ is an abelian subalgebra. Furthermore, we define
the derived bracket on the graded vector space $C^{\ast}(M, A)$ by
\begin{eqnarray*}
[[P,P']] :=(-1)^{m}[d(P), P']= (-1)^{m}[[\mu+l+r,P],P'],
\end{eqnarray*}
for any $P\in C^{m}(M, A), P'\in C^{n}(M, A)$.
The derived bracket $[[\c, \c]]$ is closed on $C^{\ast}(M, A)$, which implies that ($C^{\ast}(M, A), [[\c, \c]]$) is a graded Lie algebra.

For $T\in C^1(M, A)$, we have
\begin{eqnarray*}
&& [[T, T]](u, v)=2(Tu\c Tv-T(l(Tu)v)-T(r(Tv)u)).
\end{eqnarray*}

Thus, $T$ is a  Maurer-Cartan element  (i.e. $[[T, T]] = 0$) if and only if $T$ is a  Rota-Baxter operator on an anti-flexible algebra $(A, \c)$ with respect to the bimodule $(M, l, r)$. The proof is finished.  \hfill $\square$

Thus, Rota-Baxter operators can be characterized as Maurer-Cartan elements in a gLa. It
follows from the above theorem that if $T$ is a Rota-Baxter operator, then $d_{T}:= [[T, \c]]$
is a differential on $C^{\ast}(M, A)$ and makes the gLa $(C^{\bullet}(M, A), [[\c , \c]])$ into a dgLa.

The cohomology of the cochain complex $(C^{\bullet}(M, A), d_{T})$ is called the cohomology of the Rota-Baxter operator $T$ on $(A, \c)$ with respect to the bimodule $(M, l, r)$. We denote the corresponding cohomology groups simply by $H^{\bullet}(M, A)$.

\begin{theorem}
Let $T$ be a Rota-Baxter operator on an anti-flexible algebra $(A, \c)$ with respect to the bimodule $(M, l, r)$.  The sum $T + T'$ is a Rota-Baxter operator  if and only if $T'$ is a Maurer-Cartan element of $(C^{\ast}(M, A), [[\c, \c]], d_{T})$, that is,
\begin{eqnarray*}
[[T + T',T + T']]=0 \Leftrightarrow d_{T}T'+\frac{1}{2}[[T',T']]=0.
\end{eqnarray*}
\end{theorem}

\section{Cohomology of Rota-Baxter operators  as Hochschild cohomology }\label{sec 3}
\def\theequation{\arabic{section}. \arabic{equation}}
\setcounter{equation} {0}

In this section, we show that the cohomology of a Rota-Baxter operators on a anti-flexible algebra can also be
described as the Hochschild cohomology of a certain anti-flexible algebra with a suitable
bimodule. We also show that the cohomology of a Rota-Baxter operator on an anti-flexible algebra is related
with the cohomology of the corresponding Rota-Baxter operator on the commutator Lie algebra.

Let $T : M \rightarrow A$ be a Rota-Baxter operator on an anti-flexible algebra $(A, \c)$ with respect to the bimodule $(M, l, r)$. By
Proposition 2.12,   then  the vector space $M$ carries an anti-flexible  algebra structure with the
product
\begin{eqnarray*}
m\star_T n = r(T(n))m + l(T(m))n,~\text{for all}~ m, n\in M.
\end{eqnarray*}
\begin{lemma}
Let $T : M \rightarrow A$ be a Rota-Baxter operator on an anti-flexible algebra $(A, \c)$ with respect to the bimodule $(M, l, r)$. Define
\begin{eqnarray*}
&&l_T:M\rightarrow gl(A),~~~~l_T(m)(a)=T(m)\c a-T(r(a)m), \\
&& r_T:M\rightarrow gl(A),~~~~r_T(m)(a)=a\c T(m)-T(l(a)m), ~\text{for all}~ m\in M, a\in A.
\end{eqnarray*}
Then $l_T , r_T$ defines an $M$-bimodule structure on $(A, \c)$.
\end{lemma}
{\bf Proof.} For any $m, n \in M$ and $a \in A$, we have
\begin{eqnarray*}
&&[l_T(m\star_T n)-l_T(m)l_T(n) - r_T(m)r_T(n)+r_T(n \star_T m)](a)\\
&=&l_T(m\star_T n)(a)-l_T(m)l_T(n)(a) - r_T(m)r_T(n)(a)+r_T(n \star_T m)(a)\\
&=&(T(m)\c T(n))\c a-T(r(a)(r(T(n))m + l(T(m))n))\\
&&-T(m)\c(T(n)\c a)+T((l(T(m))(r(a)n))+r(T(n)a)m)\\
&&-(a\c T(n))\c T(m)+T(l(aT(n))m+r(T(m))l(a)n)\\
&&+a\c (T(n)\c T(m))-T(l(a)l(T(n))m-l(a)r(T(m))n)\\
&=&0.
\end{eqnarray*}
Similarly, we have
\begin{eqnarray*}
&&l_T(m)r_T(n)- r_T(n)l_T(m) - l_T(n)r_T(m)+r_T(m)l_T(n)=0.
\end{eqnarray*}
Then $l_T , r_T$ defines an $M$-bimodule structure on $(A, \c)$. Hence the proof is finished.\hfill $\square$

By Lemma 3.1 we obtain an $M$-bimodule structure on the vector space $(A, \c)$. Therefore,
we may consider the corresponding Hochschild cohomology of $M$ with coefficients in $(A, l_T, r_T)$.
More precisely, we define
\begin{eqnarray*}
C^{n}(M,A) := \text{Hom}(M^{\otimes n}, A),~\text{for all}~n \geq 0,
\end{eqnarray*}
and the differential is given by
\begin{eqnarray*}
d_H(a)(m)& =& l_T (m)(a)-r_T (m)(a) \\
&=&T(m)\c a-T(r(a)m)- a\c T(m) + T(l(a)m), ~\text{for all}~ a\in A =C^{0}(M,A),
\end{eqnarray*}
and
\begin{eqnarray*}
&&(d_Hf)(u_1,\ldots,u_{n+1}) \\
&=&T(u_1)f(u_2,\ldots, u_{n+1}) - T(r(f(u_2,\ldots, u_{n+1})u_1)\\
&&+ \sum_{i=1}^{n}\sum_{\sigma\in \Sigma_{n+1}}(-1)^{i} (-1)^{\epsilon(\sigma)}f(u_{\sigma(1)},\ldots,u_{\sigma(i-1)}, r(T(u_{\sigma(i+1)}) u_{\sigma(i)})\\
&&\hspace{8cm}+ l(T(u_{\sigma(i)}))u_{\sigma(i+1)},\ldots,u_{\sigma(n+1)})\\
&&+ (-1)^{n+1} f(u_1,\ldots,u_n)\c T(u_{n+1})-(-1)^{n+1}T(l_T(u_{n+1})f(u_1,\ldots,u_n)).
\end{eqnarray*}
We denote the group of $n$-cocycles by $Z^n(M, A)$ and the group of $n$-coboundaries
by $B^n(M, A)$. The corresponding cohomology groups are defined by $H^n(M, A) = Z^n(M, A)$\\$/B^n(M, A), n \geq 0$.

It follows from the above definition that
\begin{eqnarray*}
H^{0}(M,A)&=& \{a \in A| d_H(a)=0\}\\
&=&\{a \in A| a \c T(m)-T(m)\c a=T(l(a)m - r(a)m),~\text{for all}~m \in M\}.
\end{eqnarray*}
By \cite{DBH20}, if $a, b\in A$, define the commutator by  $[a, b]_g=a \c b-b\c a$, then it is a Lie algebra and we denote it by $(g(A), [\c , \c]_g)$.
Furthermore, it is easy to check that $H^0(M, A)$ has a Lie algebra structure induced from that of $(A, \c)$.

Note that a linear map $ f\in C^1(M, A)$ is closed if it satisfies
\begin{eqnarray*}
T(u) \c f(v) + f(u) \c T(v)- T(l_T(u)f(v)+ r_T(v)f(u))-f(l_T(u)T(v) +r_T(v)T(u))=0,
\end{eqnarray*}
for any $u, v \in M$.

For a Rota-Baxter operator $T$ on an anti-flexible algebra $(A, \c)$ with respect to the bimodule $(M, l, r)$, we get two coboundary operators $d_T =[[T, \c]]$ and $d_H$ on the same graded vector space $ C^{\bullet}(M, A) = \oplus_{n\geq 0}C^n(M, A)$.

The following proposition relates the above two coboundary operators.
\begin{proposition}
Let $T: M \rightarrow A$ be a Rota-Baxter operator on an anti-flexible algebra $(A, \c)$ with respect to the bimodule $(M, l, r)$. Then the two coboundary operators are related by
\begin{eqnarray*}
d_Tf = (-1)^n d_Hf,~\text{for all}~ \in C^{n}(M, A).
\end{eqnarray*}
\end{proposition}
{\bf Proof.}  For any $f \in C^n(M, A)$ and $u_1, \ldots, u_{n+1} \in M$, we have
\begin{eqnarray*}
&&(d_T f)(u_1,\ldots, u_{n+1}) \\
&=& [[T, f]]\\
&=&T(l(f(u_2,\ldots, u_{n+1})u_{n+1})- (-1)^{n}T(r(f(u_2,\ldots, u_{n+1})u_1) \\
&&- (-1)^{n}\{\sum_{i=1}^{n}\sum_{\sigma\in \Sigma_{n+1}}(-1)^{i-1} (-1)^{\epsilon(\sigma)}f(u_{\sigma(1)},\ldots,u_{\sigma(i-1)}, r(T(u_{\sigma(i+1)}) u_{\sigma(i)})\\
&&\hspace{7cm}+ l(T(u_{\sigma(i)}))u_{\sigma(i+1)},\ldots,u_{\sigma(n+1)})\}\\
&&+ (-1)^{n}T(u_{1})\c f(u_2,\ldots,u_{n+1}) -f(u_1,\ldots,u_n)\c T(u_{n+1})\\
&=& (-1)^{n} (d_Hf)(u_1,\ldots,u_{n+1}).
\end{eqnarray*}
This complete the proof.\hfill $\square$

\begin{definition}(\cite{G12})
Let $(g, [\c, \c]_g)$ be a Lie algebra and $\rho : g \rightarrow gl(M)$ be a representation of $g$ on a vector
space $M$. A Rota-Baxter operator  on $g$ with respect to the representation $M$ is a linear map
$T : M \rightarrow g$ satisfying
\begin{eqnarray*}
[T(m), T(n)] = T(\rho(Tm)(n)-\rho(Tn)(m)), ~\text{for all}~ m, n \in M.
\end{eqnarray*}
\end{definition}

\begin{lemma}(\cite{DBH20})
Let $(M, l, r)$ be a bimodule of an anti-flexible algebra $(A, \c)$. Then $(M, l-r)$ is a
representation of the associated Lie algebra $(g(A), [\c,\c]_g)$.
\end{lemma}

With the above notations, we have the following
\begin{proposition}
The collection of maps $S_n : \text{Hom}(A^{\otimes n}, M) \rightarrow  Hom(\wedge ^n A, M)$ defined by
\begin{eqnarray*}
S_n(f)(a_1, \ldots, a_n)=\sum_{\sigma\in \Sigma_{n}}(-1)^{\epsilon(\sigma)}f(a_{\sigma(1)},\ldots,a_{\sigma(n)})
\end{eqnarray*}
is a morphism from the Hochschild cochain complex of $A$ with coefficients in the bimodule
$M$ to the Chevalley-Eilenberg complex of the commutator Lie algebra $(g(A), [\c,\c]_g)$ with coefficients
in the representation $(M, l-r)$.
\end{proposition}

\begin{proposition}
Let $T: M \rightarrow A$ be a Rota-Baxter operator on an anti-flexible algebra $(A, \c)$ with respect to the bimodule $(M, l, r)$. Then $T$ is also  a Rota-Baxter operator on the commutator Lie algebra $(g(A), [\c , \c]_g)$ with respect to the representation $(M, \rho)$.
\end{proposition}
{\bf Proof.} For any $m, n\in M$, we have
\begin{eqnarray*}
[T(m), T(n)] &=& T(m) \c T(n)- T(n) \c T(m)\\
&=& T(r(T(n))m+l(T(m))n)- T(r(T(m))n + l(T(n))m)\\
&=& T((l-r)(T(m)n-(l-r)(T(n)m)\\
&=& T(\rho(T(m)n-\rho(T(n)m).
\end{eqnarray*}
This proves the proposition. \hfill $\square$
\section{Infinitesimal deformations of bimodules over anti-flexible algebras} \label{sec 4}
\def\theequation{\arabic{section}. \arabic{equation}}
\setcounter{equation} {0}

Let $(A, \c)$ be an anti-flexible algebra and $(M, l, r)$ a  bimodule. Let $\omega: \otimes^2A \rightarrow A$, $\phi: A\rightarrow gl(M)$ and $\psi: A\rightarrow gl(M)$   be  linear maps.  Consider a $t$-parametrized family of multiplication operations and linear maps:
\begin{eqnarray*}
a\c_t b= a\c b+\omega(a, b),~~l^{t}(a)=l(a)+t\phi(a),~~~r^{t}(a)=r(a)+t\psi(a),~\text{for all}~ a, b\in A.
\end{eqnarray*}
If $(A, \c_t)$ are anti-flexible algebras and  $(M, l^t, r^t)$ are   bimodules for all $t\in \mathbb{K}$,    we say that
$(\omega, \phi, \psi)$ generates an infinitesimal deformation of the $A$-bimodule $M$.

Let $\mu_t$ denote the anti-flexible algebra structure $(A, \c_t)$. By Proposition 2.16, the bimodule $(M, l^t, r^t)$
over the anti-flexible algebra $(A, \c_t)$ is an infinitesimal deformation of the $A$-bimodule $M$ if and only if
\begin{eqnarray*}
(\mu_t+l^t+r^t)\overline{\circ}(\mu_t+l^t+r^t)=0,
\end{eqnarray*}
which is equivalent to
\begin{eqnarray*}
&&(\omega+\phi+\psi)\overline{\circ}(\omega+\phi+\psi)=0,\\
&& (\mu_t+l^t+r^t)\overline{\circ}(\omega+\phi+\psi)=0.
\end{eqnarray*}
\begin{definition}
Let the bimodules $(M, l^t, r^t)$ and $(M, l'^t, r'^t)$ be two infinitesimal deformations of an $A$-bimodule $M$ over anti-flexible algebras $(A, \c_t)$ and
$(A, \c'_t)$ respectively. We call them equivalent
if there exists $N \in gl(g)$ and $S \in gl(M)$ such that $(Id_A + tN, Id_M + tS )$ is a homomorphism from the
bimodule $(M, l'^t, r'^t)$ to the bimodule $(M, l^t, r^t)$.
\end{definition}
By direct calculations, the bimodule $(M, l^t, r^t)$ over the anti-flexible algebra $(A, \c_t)$ and the bimodule
$(M, l'^t, r'^t)$ over the anti-flexible algebra $(A, \c'_t)$ are equivalent deformations if and only if
\begin{eqnarray*}
&& ( \omega+ \phi +\psi)(a + m, b + n) -( \omega'+ \phi' +\psi')(a + m, b + n)  = d(N + S)(a + m, b + n),\\
&&( \omega'+ \phi' +\psi')(N(a) +S (m), N(b) + S (n))= 0,
\end{eqnarray*}
and
\begin{eqnarray*}
&&(N + S)( \omega+ \phi +\psi)(a + m, b + n) = ( \omega'+ \phi' +\psi')(a + m, N(b) + S(n))\\
&&+( \omega'+ \phi' +\psi')(N(a) + S(m), b + n)+ ( \mu + l + r)(N(a) + S (m), N(b) + S (n)).
\end{eqnarray*}
Summarizing the above discussion, we have the  following conclusion:
\begin{theorem}
Let the bimodule $(M, l^t, r^t)$  over the anti-flexible algebra $(A, \c_t)$ be an infinitesimal
deformation of an $A$-bimodule $M$ generated by $(\omega, \phi, \psi)$. Then $\omega+\phi+\psi \in C^2(M, A)$ is closed, i.e.
$d(\omega, \phi, \psi)= 0$. Furthermore, if two infinitesimal deformations $(M, l^t, r^t)$ and $(M, l'^t, r'^t)$ over
anti-flexible algebra $(A, \c_t)$ and $(A, \c'_t)$ generated by $(\omega, \phi, \psi)$ and $(\omega', \phi', \psi')$ respectively are equivalent,
then $\omega+\phi+\psi $ and $\omega'+\phi'+\psi' $ are in the same cohomology class in $H^2(M, A)$.
\end{theorem}
\begin{definition}
An infinitesimal deformation of an $A$-bimodule $M$ is said to be trivial if it is equivalent to the
$A$-bimodule $M$.
\end{definition}
One can deduce that the bimodule $(M, l^t, r^t)$ over the anti-flexible algebra $(A, \c_t)$ is a trivial deformation
if and only if for all $a, b \in A, m, n \in M$, we have
\begin{eqnarray*}
&&( \omega+\phi+\psi)(a+m, b+n)=d(N+S)(a+m, b+n),\\
&&(N+S)( \omega+\phi +\psi)(a+ m, b+n) =(\mu+l+r)(N(a)+S (m), N(b)+S (n)).
\end{eqnarray*}
Equivalently, we have
\begin{eqnarray}
&&\omega(a, b) = N(a)\c b + a \c N(b)- N(a \c b),\\
&& N\omega(a, b) = N(a) \c N(b),\\
&&\phi(a) = l(N(a)) + l(a)\circ  S - S \circ l(a),\\
 &&l(N(a)) \circ S = S \circ \phi(a),\\
&&\psi(a) = r(N(a)) + r(a) \circ S - S \circ r(a),\\
&&r(N(a)) \circ S = S \circ \psi(a).
\end{eqnarray}
It follows from Eqs.(4.1) and (4.2) that $N$ must be a Nijenhuis operator on the anti-flexible algebra $(A, \c)$. It
follows from Eqs.(4.3) and (4.4) that $N$ and $S$ should satisfy the condition:
\begin{eqnarray}
 &&l(N(a))S(m) = S( l(N(a))(m) + l(a)(S(m)) - S(l(a)m)), ~\forall a\in A, m\in M.
\end{eqnarray}
It follows from Eqs. (4.5) and (4.6) that N and S should also satisfy the condition:
\begin{eqnarray}
&&r(N(a)) S(m) = S(r(N(a))m + r(a)S(m) - S (r(a)m)), ~~\forall a\in A, m\in M.
\end{eqnarray}
\begin{theorem}
Let $(M, l, r)$ be a bimodule over an anti-flexible algebra $(A, \c)$, $N \in gl(A)$ and $S \in  gl(M)$.
If $N$ is a Nijenhuis operator on the anti-flexible algebra $(A, \c)$ and if $S$ satisfies Eqs. (4.7) and (4.8), then a
deformation of the $A$-bimodule $M$ can be obtained by putting
\begin{eqnarray*}
\omega(a, b) = N(a)\c  b + a \c N(b)- N(a\c b),\\
\phi(a) = l(N(a)) + l(a) \circ S - S \circ l(a),\\
\psi(a) = r(N(a)) + r(a) \circ S - S \circ r(a),
\end{eqnarray*}
for any $a, b\in A$. Furthermore, this deformation is trivial.
\end{theorem}
Note that the conditions that $N$ is a Nijenhuis operator and $S$ satisfies Eqs. (4.7) and (4.8), can be expressed
simply by the following result.
\begin{proposition}
Let $(M, l, r)$ be a bimodule over an anti-flexible algebra $(A, \c)$. Then $N$ is a Nijenhuis
operator on the anti-flexible algebra $(A, \c)$ and $S$ satisfies satisfies Eqs. (4.7) and (4.8), if and only if $N + S$ is a
Nijenhuis operator on the semidirect product anti-flexible algebra $A\oplus M$.
\end{proposition}

\begin{definition}
Let $(M, l, r)$ be a bimodule over an anti-flexible algebra $(A, \c)$. A pair $(N, S )$, where
$N \in gl(A)$ and $S \in  gl(M)$, is called a Nijenhuis structure on an $A$-bimodule $M$ if $N$ and $S^{\ast}$ generate a trivial
infinitesimal deformation of the dual $A$-module $M^{\ast}$.
\end{definition}
Note that the condition of the above definition is equivalent to the fact that $N$ is a Nijenhuis tensor on
$A$ and
\begin{eqnarray*}
&&l(N(a))S(m) = S (l(N(a))m) + l(a)S^2(m)-S (l(x)S (m)),\\
&&r(N(a))S(m)= S (r(N(a))m) + r(a)S^2(m)- S (r(x)S (m)).~\text{for all}~ a\in A, m\in M.
\end{eqnarray*}
\begin{example}
Let $N: A\rightarrow A$ be a Nijenhuis operator on the anti-flexible algebra $(A, \c)$. Then $(N, N^{\ast})$ is a Nijenhuis structure on
the coadjoint module $A^{\ast}$.
\end{example}
\begin{cor}
Let $(N, S)$  be a Nijenhuis structure on a $A$-bimodule $M$, then the pairs $(N^{i}, S^{i})$  are  Nijenhuis structures on an $A$-bimodule $M$.
\end{cor}

\section{$\mathcal{ON}$-structures on bimodules over anti-flexible algebras and compatible Rota-Baxter operators}\label{sec 5}
\def\theequation{\arabic{section}. \arabic{equation}}
\setcounter{equation} {0}
In this final section, we show how compatible Rota-Baxter operators and $\mathcal{ON}$-structures are related.

Let $T : M \rightarrow A$ be a Rota-Baxter operator on an anti-flexible algebra $(A, \c)$ with respect to the bimodule $(M, l, r)$. By
Proposition 2.12,   then  the vector space $M$ carries an anti-flexible  algebra structure with the
product
\begin{eqnarray*}
m\star_T n = l(T(m))n +r(T(n))m,  ~\text{for all}~ m, n\in M.
\end{eqnarray*}
We define the multiplication $\star^S_{T}: M\otimes M\rightarrow M$ to be the deformed multiplication of $\star_T$ by $S$ , i.e.
\begin{eqnarray*}
m\star^S_{T} n=S (m)\star_T n + m \star_T S (n)-S (m \star_T n).
\end{eqnarray*}
\begin{definition}
Let $T : M \rightarrow A$ be a Rota-Baxter operator and $(N, S )$ a Nijenhuis structure on an $A$-bimodule $M$. The triple $(T, N, S )$ is called an $\mathcal{ON}$-structure on an $A$-bimodule $M$ if $T$ and $(N, S)$
satisfy the following conditions
\begin{eqnarray*}
N \circ T = T \circ S,\\
m\star_{N \circ T} n=m\star^S_{T} n,~\text{for all}~ m, n\in M.
\end{eqnarray*}
\end{definition}
Define two linear maps $\widetilde{l}, \widetilde{r}: A\rightarrow gl(M)$ as follows:
\begin{eqnarray*}
&&\widetilde{l}(a):=l(N(a))-l(a)\circ S+S\circ l(a),\\
&&\widetilde{r}(a):=r(N(a))-r(a)\circ S+S\circ r(a), ~\text{for all}~ a\in A.
\end{eqnarray*}
Then it is easy to check that $(M, \widetilde{l}, \widetilde{r})$ is a bimodule of $(A, \c)$. Furthermore, we have an anti-flexible  algebra structure with the
product
\begin{eqnarray*}
m\widetilde{\star}_T n = \widetilde{l}(T(m))n+\widetilde{r}(T(n))m ,~\text{for all}~m, n\in M.
\end{eqnarray*}
Direct calculation, for any $m, n\in M$,  we have $m\widetilde{\star}_T n+m\star^S_{T} n=2(m\star_{N \circ T} n)$. Then we have the following lemma.
\begin{lemma}
Let $(T, N, S)$ be an an $\mathcal{ON}$-structure on an anti-flexible algebra $(A, \c)$ with respect to the bimodule $(M, l, r)$. Then we have
\begin{eqnarray*}
m\star^S_{T} n=m\widetilde{\star}_{T} n.
\end{eqnarray*}
\end{lemma}
\begin{definition}
Two Rota-Baxter operators $T_1, T_2: M \rightarrow A$ on an anti-flexible algebra $(A, \c)$ with respect to the bimodule $(M, l, r)$ are said to be compatible if their sum $T_1 + T_2: M \rightarrow A$ is also a Rota-Baxter operator.
\end{definition}
\begin{proposition}
Let $T_1, T_2: M \rightarrow A$ be two Rota-Baxter operators on an anti-flexible algebra $(A, \c)$ with respect to the bimodule $(M, l, r)$ . If $T_1, T_2$ are compatible and $T_2$ is invertible then $N = T_1 \circ T_2^{-1}: A\rightarrow A$ is a Nijenhuis operator on the anti-flexible algebra $(A, \c)$.
Conversely, if $T_1, T_2$ are both invertible and $N$ is a Nijenhuis tensor then $T_1, T_2$ are compatible.
\end{proposition}
{\bf Proof.} Let $T_1, T_2$ be compatible and $T_2$ invertible. For any $a, b \in A$, there exists elements $m, n \in M$
such that $T_2(m) = a$ and $T_2(n) = b$. Then
\begin{eqnarray*}
&& Na\c Nb - N(Na\c b + a\c Nb) + N^2(a\c b)\\
&=& NT_2(m)\c NT_2(n) - N(NT_2(m)\c T_2(n) + T_2(m)\c NT_2(n)) + N^2(T_2(m)\c T_2(n))\\
&=& T_1(m)\c T_1(n) - N(T_1(m)\c T_2(n) + T_2(m)\c T_1(n)) + N^2(T_2(m)\c T_2(n))\\
&=& T_1(l(T_1(m))n + r(T_1(n))m)-N(T_1(l(T_2(m))n + r(T_2(n))m)\\
&&+T_2(l(T_1(m))n + r(T_1(n))m))+ N^{2}(T_2(l(T_2(m))n + r(T_2(n))m))\\
&=& 0.
\end{eqnarray*}
Conversely, if $N$ is a Nijenhuis tensor then for all $m, n \in M$,
\begin{eqnarray*}
NT_2(m)\c NT_2(n) = N(NT_2(m)\c T_2(n) + T_2(m)\c NT_2(n)) - N^2(T_2(m)\c T_2(n)).
\end{eqnarray*}
This implies that
\begin{eqnarray*}
&& T_1(l(T_1(m))n + r(T_1(n))m)\\
&=&N(T_1(m)\c T_2(n)+T_2(m)\c T_1(n))- NT_1(l(T_2(m))n + r(T_2(n))m).
\end{eqnarray*}
Since $N$ is invertible, we may apply $N^{-1}$ to both sides to get the above identity. Hence $T_1$ and $T_2$ are
compatible.  \hfill $\square$
\begin{theorem}
Let $(T, N, S)$ be an $\mathcal{ON}$-structure on an anti-flexible algebra $(A, \c)$ with respect to the bimodule $(M, l, r)$. Then\\
(i) $T$ is a Rota-Baxter operator on  the deformed anti-flexible algebra $(A, \c_N)$ with respect to the bimodule $(M, \widetilde{l}, \widetilde{r})$, \\
(ii) $N \circ T$ is a Rota-Baxter operator on an anti-flexible algebra $(A, \c)$ with respect to the bimodule $(M, l, r)$.
\end{theorem}
{\bf Proof.} (i) For any $m, n \in  M$, we have
\begin{eqnarray*}
T(m\widetilde{\star}_Tn)&=& T(m\star^{S}_Tn)\\
&=& T(S (m)\star_T n + m \star_T S (n)-S (m \star_T n))\\
&=& T\circ S (m)\c_N T(n) + T(m) \c_N T\circ S (n)-T\circ S (m \star_T n)\\
&=&  N\circ T (m)\c_N T(n) + T(m) \c_N N\circ T (n)-N (T(m) \c_N  T(n))\\
&=& T(m) \c_N T(n).
\end{eqnarray*}
Then $T$ is a Rota-Baxter operator on  the deformed anti-flexible algebra $(A, \c_N)$ with respect to the bimodule $(M, \widetilde{l}, \widetilde{r})$.\\
(ii) By the fact that $N$ is a Nijenhuis tensor, we have
\begin{eqnarray*}
N\circ T(m \star_{N \circ T} n) = N\circ T(m\star^{S}_Tn) = N(T(m)\c_N T(n)) = N\circ T(m) \c N\circ T(n).
\end{eqnarray*}
Hence  $N \circ T$ is a Rota-Baxter operator on an anti-flexible algebra $(A, \c)$ with respect to the bimodule $(M, l, r)$. \hfill $\square$
\begin{proposition}
Let $(T, N, S)$ be an $\mathcal{ON}$-structure on an anti-flexible algebra $(A, \c)$ with respect to the bimodule $(M, l, r)$. Then $T$ and $N \circ T$ are compatible Rota-Baxter operators.
\end{proposition}
{\bf Proof.} For any $m, n \in M$, we have
\begin{eqnarray*}
&&  m \star_{T+N \circ T} n=m \star_{T} n+m \star_{N \circ T} n=m \star_{T} n+m \star^{S}_{T} n.
\end{eqnarray*}
Furthermore, we have
\begin{eqnarray*}
&&(T+N \circ T)(m \star_{T+N \circ T} n)\\
&=& T(m \star_{T} n)+T(m \star^{S}_{T} n)+(N\circ T)(m \star_{T} n)+(N\circ T)(m \star^{S}_{T} n)\\
&=&T(m \star_{T} n)+T(S (m)\star_T n + m \star_T S (n)-S (m \star_T n))\\
&&+(N\circ T)(m \star_{T} n)+(N\circ T)(m \star_{N\circ T} n)\\
&=&T(m) \c T(n)+(T\circ S) (m)\c T(n) + T(m) \c (T\circ S) (n)+(N\circ T)(m) \c  (N\circ T)(n)\\
&=&T(m) \c T(n)+(N\circ T) (m)\c T(n) + T(m) \c (N\circ T) (n)+(N\circ T)(m) \c  (N\circ T)(n)\\
&=&(T+N \circ T)(m)\c (T+N \circ T)(n).
\end{eqnarray*}
Then $T$ and $N \circ T$ are compatible Rota-Baxter operators.  \hfill $\square$

In the next proposition, we construct an $\mathcal{ON}$-structure from compatible Rota-Baxter operators.
\begin{proposition}
Let $T_1, T_2 : M \rightarrow A$ be two compatible Rota-Baxter operators on an anti-flexible algebra $(A, \c)$ with respect to the bimodule $(M, l, r)$. If $T_2$ is invertible then $(T_2, N = T_1 \circ T^{-1}_2, S = T^{-1}_2\circ T_1)$ is an $\mathcal{ON}$-structure.
\end{proposition}
{\bf Proof.} Since $T_1, T_2$ are compatible Rota-Baxter operators,       we have
\begin{eqnarray*}
T_1(m)\c T_2(n)+T_2(m)\c T_1(n)=T_1(l(T_2(m))n +r(T_2(n))m)+T_2(l(T_1(m))n +r(T_1(n))m).
\end{eqnarray*}
By replacing $T_1$ with $T_2 \circ S$ in above equation, we have
\begin{eqnarray}
&&(T_2 \circ S)(m)\c T_2(n)+T_2(m)\c (T_2 \circ S)(n)\nonumber\\
&=&(T_2 \circ S)(l(T_2(m))n +r(T_2(n))m)+T_2(l((T_2 \circ S)(m))n +r((T_2 \circ S)(n))m).~~~~
\end{eqnarray}
On the other hand, $T_2$ is a Rota-Baxter operator implies that
\begin{eqnarray}
&&(T_2 \circ S)(m)\c T_2(n)+T_2(m)\c (T_2 \circ S)(n)\nonumber\\
&=&T_2( l(T_2(S(m)))n +r(T_2(n))S(m)+l(T_2(m))S(n) +r(T_2(S(n)))m).
\end{eqnarray}
From Eqs.(5.1) and (5.2) and using the fact that $T_2$ is invertible, we get
\begin{eqnarray*}
S(l(T_2(m))n +r(T_2(n))m)=l(T_2(m))S(n)+r(T_2(n))S(m).
\end{eqnarray*}
By replacing $n$ by $S(n)$, we have
\begin{eqnarray}
S(l(T_2(m))S(n) +r((T_2\circ S)(n))m)=l(T_2(m))S^{2}(n)+r(T_2\circ S(n))S(m).
\end{eqnarray}
As $T_1 = T_2 \circ S$ and $T_2$ are Rota-Baxter operators,
\begin{eqnarray*}
T_2(m\star_{T_2\circ S}n)=(T_2\circ S)(m)\c (T_2\circ S)(n)=T_2(S(m)\star_{T_2}S(n)).
\end{eqnarray*}
The invertibility of $T_2$ implies that
\begin{eqnarray}
S(l((T_2\circ S)(m))n +r((T_2\circ S)(n))m)=l((T_2\circ S)(m))S(n)+r((T_2\circ S)(n))S(m).
\end{eqnarray}
From Eqs.(5.3) and (5.4) and using the fact that $T_2$ is invertible, we get
\begin{eqnarray*}
l((T_2\circ S)(m))S(n)=l(T_2(m))S^{2}(n)+S(l((T_2\circ S)(m))n-S(l(T_2(m))S(n).
\end{eqnarray*}
Substitute $a = T_2(m)$, using $T_2\circ S = N\circ T_2$ and the invertibility of $T_2$,
\begin{eqnarray*}
l(N(a))S(n)=l(a)S^{2}(n)+S(l(N(a))n-S((l(a)S(n)).
\end{eqnarray*}
Hence the identity Eq. (4.7) follows. Similarly, Eq. (4.8) holds.  Thus, the pair $(N, S)$ is a Nijenhuis structure on an anti-flexible algebra $(A, \c)$ with respect to the bimodule $(M, l, r)$.

Next, observe that $N \circ T_2=T_2\circ S=T_1$. Moreover,
\begin{eqnarray*}
&&m\star^{S}_{T_2}n-m\star_{T_2\circ S}n\\
&=& l(T_2(m))S(n)+r(T_2(n))S(m)-S(l(T_2(m))n+r(T_2(n))m)\\
&=&0,
\end{eqnarray*}
which implies that $m\star^{S}_{T_2}n=m\star_{T_2\circ S}n$. Therefore, $(T_2, N = T_1 \circ T^{-1}_2, S = T^{-1}_2\circ T_1)$ is an $\mathcal{ON}$-structure.   \hfill $\square$
\begin{center}
 {\bf ACKNOWLEDGEMENT}
 \end{center}

  The paper is  supported by the NSF of China (No. 11761017), Guizhou Provincial  Science and Technology  Foundation (No. [2020]1Y005).

\renewcommand{\refname}{REFERENCES}

\end{document}